\documentstyle[amstex,amscd,ctagsplt,twoside]{amsart}

\catcode`@=11
\atdef@ I#1I#2I{\CD@check{I..I..I}{\llap{$\m@th\vcenter{\hbox
  {$\scriptstyle#1$}}$}
  \rlap{$\m@th\vcenter{\hbox{$\scriptstyle#2$}}$}&&}}
\atdef@ E#1E#2E{\ampersand@
  \ifCD@ \global\bigaw@\minCDarrowwidth \else \global\bigaw@\minaw@ \fi
  \setboxz@h{$\m@th\scriptstyle\;\;{#1}\;$}%
  \ifdim\wdz@>\bigaw@ \global\bigaw@\wdz@ \fi
  \@ifnotempty{#2}{\setbox@ne\hbox{$\m@th\scriptstyle\;\;{#2}\;$}%
    \ifdim\wd@ne>\bigaw@ \global\bigaw@\wd@ne \fi}%
  \ifCD@\enskip\fi
    \mathrel{\mathop{\hbox to\bigaw@{}}%
      \limits^{#1}\@ifnotempty{#2}{_{#2}}}%
  \ifCD@\enskip\fi \ampersand@}
\def\mathunderline#1{\relax\protect\@@underline{#1}}
\catcode`@=\active

\newcommand{\site}{{\mathcal E}}

\newcommand{\Hom}{\operatorname{Hom}}
\newcommand{\HOM}{\operatorname{\mathcal H{\frak om}}}

\newcommand{\Ext}{\operatorname{Ext}}

\newcommand{\Pic}{\operatorname{Pic}}
\newcommand{\NS}{\operatorname{NS}}
\newcommand{\Alb}{\operatorname{A}}

\newcommand{\Spec}{\operatorname{Spec}}

\newcommand{\coker}{\operatorname{coker}}

\newcommand{\id}{\operatorname{id}}
\newcommand{\Br}{\operatorname{Br}}
\newcommand{\cd}{\operatorname{cd}}
\newcommand{\Ab}{\operatorname{Ab}}
\newcommand{\Mod}{\operatorname{Mod}}
\newcommand{\CH}{\operatorname{CH}}
\renewcommand{\lim}{\operatorname{lim}}
\newcommand{\colim}{\operatornamewithlimits{colim}}
\newcommand{\rk}{\operatorname{rank}}
\newcommand{\corank}{\operatorname{corank}}

\newtheorem{lemma}{Lemma}[section]
\newtheorem{prop}[lemma]{Proposition}
\newtheorem{theorem}[lemma]{Theorem}

\newtheorem{cor}[lemma]{Corollary}
\newtheorem{conj}[lemma]{Conjecture}
\newcommand{\proof}{\noindent{\it Proof.\ }}
\newcommand{\proofend}{\hfill $\Box$ \\}

\newcommand{\example}{\medskip\noindent{\bf Example.\ }}
\newcommand{\f}{{\cal F}}
\newcommand{\g}{{\cal G}}
\newcommand{\N}{{\Bbb N}}

\newcommand{\Z}{{{\Bbb Z}}}
\newcommand{\Q}{{{\Bbb Q}}}

\newcommand{\F}{{{\Bbb F}}}
\newcommand{\fun}{{{\Xi}}}
\newcommand{\T}{{\cal T}}
\newcommand{\M}{{\cal M}}
\renewcommand{\L}{{\cal L}}

\newcommand{\et}{{\text{\rm\'et}}}
\newcommand{\Zar}{{\text{\rm Zar}}}

\renewcommand{\mod}{\operatorname{mod}}
\begin{document}

\title{Weil-\'etale cohomology over finite fields}

\author{Thomas Geisser*}

\thanks{${}^{*\phantom{*}}$\ Supported in part by JSPS,
the Alfred P. Sloan Foundation, and NSF
Grant No. 0070850}

\address{University of Southern California\\
Department of Mathematics, DRB \\
1042 W. 36th Place \\
Los Angeles, CA 90089}


\begin{abstract}
We calculate the derived functors $R\gamma_*$ for the base change
$\gamma$
from the Weil-\'etale site to the \'etale site for a variety over
a finite field. For smooth and proper varieties, we apply this to
express Tate's conjecture and Lichtenbaum's conjecture on special
values of $\zeta$-functions in terms of Weil-\'etale cohomology of the
motivic complex $\Z(n)$.
\end{abstract}

\maketitle

\section{Introduction}
In \cite{licht}, Lichtenbaum defined Weil-\'etale cohomology
groups of varieties over finite fields in order to produce finitely 
generated cohomology groups which are related to special values
of zeta functions. He gave several examples where these groups
were indeed finitely generated.
The purpose of this paper is to eludicate the precise relationship between
Weil-\'etale cohomology groups and \'etale cohomology groups. This is 
applied to give necessary and sufficient conditions for the Weil-\'etale 
cohomology groups to be finitely generated, and to be related to special
values of zeta functions.

Recall that an \'etale sheaf on a variety $X$ over a finite field $\F_q$ 
corresponds to a sheaf on $\bar X=X\times_{\F_q}\bar \F_q$, together with a 
continuous action of the Galois group $\hat G=Gal(\bar\F_q:\F_q)$. 
In the Weil-\'etale topology, the role of the Galois group
is replaced by the Weil group $G$, which is the subgroup of $\hat G$
generated by the Frobenius operator $\varphi$:
A Weil-\'etale sheaf is an \'etale sheaf on $\bar X$, together with an action 
of $G$. If we denote the category of Weil-\'etale sheaves by
$\T_G$ and the category of \'etale sheaves by $\T_{\hat G}$, then
there is a morphism of topoi $\gamma:\T_G \to \T_{\hat G}$. The
functor $\gamma^*$ is the restriction functor, and for
$U$ \'etale over $\bar X$, $\gamma_*\f(U)=\colim_H\f(U)^H$, where $H$ runs 
through sufficiently small subgroups of $G$.

We give an explicit description of the total derived functor $R\gamma_*$
and derive formulas for $\gamma_*\f$, $R^1\gamma_*\f$; for $i>1$, 
$R^i\gamma_*\f=0$. If $\f=\gamma^*\g$ is the restriction of an \'etale
sheaf, then the formula can be simplified to the following
projection formula: 

\begin{theorem}
For every complex $\g^\cdot $ of \'etale sheaves,
there is a quasi-isomorphism of complexes of \'etale sheaves
$$R\gamma_*\Z\otimes^L \g^\cdot \cong  R\gamma_*(\gamma^*\g^\cdot).$$
\end{theorem}

This raises the question of calculating $R\gamma_*\Z$. We show
that $\gamma_*\Z\cong \Z$, $R^1\gamma_*\Z\cong \Q$, which 
gives the distinguished triangle
\begin{equation}
\label{trt} \g^\cdot \to R\gamma_*\gamma^*\g^\cdot \to
\g^\cdot\otimes\Q[-1] @>\delta>> \g^\cdot[1],
\end{equation}
and implies that for a complex with torsion cohomology sheaves 
$\g^\cdot\cong R\gamma_*\gamma^*\g^\cdot$. We show that $R\gamma_*\Z$
is quasi-isomorphic to a complex considered by Kahn in \cite{kahn},
hence the map $\delta$ is induced by the composition
$$ \Q[-1]\to \Q/\Z[-1]@>\cup \epsilon>> \Q/\Z[0]@>\beta >> \Z[1],$$
with $\beta$ the Bockstein-homomorphism and 
$\epsilon\in \Ext_{\hat G}^1(\Q/\Z,\Q/\Z)$ the class of the
$\hat G$-module $\Q/\Z\oplus \Q/\Z$ with action $g\cdot(a,b)=(a+gb,b)$.
In particular, $\delta=0$ for a complex with $\Q$-vector spaces
as cohomology sheaves, hence the sequence \eqref{trt} splits and
$R\gamma_*\gamma^*\g^\cdot \cong \g^\cdot\oplus \g^\cdot[-1]$.
We show that under the latter isomorphism, the cup-product with a 
generator $e\in H^1_W(\F_q,\Z)$ is given by multiplication with the 
matrix ${0\,0\choose 1\,0}$.

In the second half of the paper, we give applications of the above 
calculations to the Weil-\'etale hypercohomology groups $H^i_W(X,\Z(n))$ 
of the motivic complex $\gamma^*\Z(n)$. We assume that $X$ is smooth over 
$\F_q$, because Weil-\'etale cohomology groups for singular schemes are 
not well-behaved. (We discuss in a forthcoming paper how to refine the 
Weil-\'etale topology to get reasonable cohomology groups for singular 
schemes). The general results above specialize to this situation, and
we show that if $X$ is of dimension $d$, then
$H^i_W(X,\Z(n))=0$ for $i>\max\{2d+1,n+d+1\}$. If
$X$ is connected and proper, then there is an isomorphism 
$H^{2d+1}_W(X,\Z(d))@>deg >> \Z$, and the composition 
$H^{2d}_W(X,\Z(d)) @>deg(-\cup e)>> \Z $ is surjective.

Lichtenbaum expected statement $L(X,n)$: If $X$ is smooth and proper,
then the cohomology groups 
$H^i_W(X,\Z(n))$ are finitely generated for all $i$. 
On the other hand, a conjecture
of Kahn \cite{kahn} can be reformulated with the above results into
statement $K(X,n)$: If $X$ is smooth and proper, then 
Weil-\'etale motivic cohomology is an integral model for $l$-adic cohomology, 
i.e. for every prime $l$ (including $p$),
$$ H^i_W(X,\Z(n))\otimes \Z_l\cong H^i_{cont}(X,\Z_l(n)).$$
Statements $K(X,n)$ and $L(X,n)$ are related to 
the conjunction $T(X,n)$ of Tate's conjecture on the surjectivity
of the cycle map 
$CH^n(X)\otimes\Q_l\to H^{2n}_{cont}(\bar X,\Q_l(n))^{\hat G}$, 
and semi-simplicity of $H^{2n}_{cont}(\bar X,\Q_l(n))$
at the eigenvalue $1$, together with Beilinson's conjecture that
rational and numerical equivalence on $X$ agree up to torsion in codimension 
$n$ (see also \cite{kahnneu}): 

\begin{theorem} Let $X$ be a smooth projective variety over $\F_q$,
and $n$ an integer. Then
$$ K(X,n)+K(X,d-n) \Rightarrow L(X,n)\Rightarrow K(X,n)
\Rightarrow T(X,n).$$
Conversely, if $T(X,n)$ holds for {\rm all} smooth and projective
varieties over $\F_q$ and {\rm all} $n$, then $K(X,n)$ and
$L(X,n)$ hold for all $X$ and $n$.
\end{theorem}

Finally we reinterpret a result of Milne \cite{values} to show
that Weil-\'etale motivic cohomology can
be used to give formulas for special values of $\zeta$-functions
of varieties over finite fields, as anticipated by Lichtenbaum. 
For a complex with finitely many finite cohomology groups, define 
$ \chi(C^\cdot) := \prod_i |H^i(C^\cdot)|^{(-1)^i}$ and let
$$ \chi(X,{\cal O}_X,n)=\sum_{i\leq n,j\leq d}(-1)^{i+j}(n-i)
\dim H^j(X,\Omega^i).$$ Since $e^2=0$, the
groups $H^*_W(X,\Z(n))$ form a complex with differential $e$.

\begin{theorem}
Let $X$ be a smooth projective variety such that $K(X,n)$
holds. Then the order $\rho_n$ of the pole of $\zeta(X,s)$
at $s=n$ is $\rk H^{2n}_W(X,\Z(n))$, and
\begin{equation*}
\zeta(X,s)=\pm
(1-q^{n-s})^{-\rho_n}\cdot \chi(H^*_W(X,\Z(n)),e)\cdot
q^{\chi(X,{\cal O}_X,n)} \quad {\text as}\; s\to n.
\end{equation*}
If furthermore $K(X,d-n)$ holds, then
\begin{equation*}
\chi(H^*_W(X,\Z(n)),e)= \prod_i
|H^i_W(X,\Z(n))_{tor}|^{(-1)^i}\cdot R^{-1},
\end{equation*}
where $R$ is the determinant of the pairing
$$ H^{2n}_W(X,\Z(n)) \times H^{2(d-n)}_W(X,\Z(d-n))
\to H^{2d}_W(X,\Z(d)) \to \Z .$$
\end{theorem}

To give explicit evidence, we show that
$K(X,0)$ holds, and
that the surjectivity of the cycle map
$\Pic X \otimes\Q_l\to H^2_{cont}(X,\Q_l(1))$
implies $K(X,1)$. In particular,
$K(X,1)$ holds for Hilbert modular
surfaces, Picard modular surfaces, Siegel modular threefolds, and
in characteristic at least $5$ for supersingular and elliptic K3
surfaces. Using the method of Soul\'e \cite{soule}, we also show
that
$K(X,n)$ holds for a smooth projective variety
$X$ of dimension $d$, which can be constructed out of products of
smooth projective curves by union, base extension and blow-ups,
and for $n\leq 1$ or $n\geq d-1$. This applies to abelian
varieties, unirational varieties of dimension at most $3$, and to
Fermat hypersurfaces. In \cite{kahnneu}, Kahn shows that
conjecture $K(X,n)$ is true for arbitrary $n$ if $X$ is of abelian
type and satisfies Tate's conjecture. This applies in particular
to the product of elliptic curves.

This paper is based on ideas of Lichtenbaum \cite{licht}
and Kahn \cite{kahn}. We wish to thank B. Kahn, S. Lichtenbaum and
T. Saito for several helpful comments.
The paper was written while the author was visiting the University
of Tokyo, which provided excellent working conditions.

\section{Profinite completion}
We fix a finite field $\F_q$, let $\bar\F_q$ be the algebraic closure 
of $\F_q$, and $\varphi$ the
arithmetic Frobenius endomorphism $x\to x^q$ of $\bar \F_q$ over
$\F_q$. The Galois group $\hat G$ of $\bar\F_q/\F_q$ is isomorphic
to the profinite completion $\lim_m\Z/m$ of $\Z$, and we let $G$
be the subgroup of $\hat G$ generated by
$\varphi$. Of course, $G$ is isomorphic to $\Z$, but we want to
avoid confusing $G$-modules and abelian groups. The fixed
field of $mG$ and of $m\hat G$ is $\F_{q^m}$. 

Let $\site$ be a full subcategory of the category of separated schemes of 
finite type over $\F_q$, which contains with every scheme $X$
also every scheme $U$ which is \'etale and of finite type over $X$. 
Our main examples
will be the category of separated schemes of finite type over $\F_q$,
the category of smooth schemes of finite type over $\F_q$,
and the small \'etale site of a scheme $X$ separated and
of finite type over $\F_q$.
Let $\bar\site$ be the full subcategory of separated schemes of finite type 
over 
$\bar \F_q$ which are connected components of the base-change of a scheme 
in $\site$; note that every scheme of finite type over $\bar \F_q$ is the
base-change of a scheme over some $\F_{q^r}$.
We equip $\site$ and $\bar\site$ with the \'etale topology, although all 
arguments below hold for any Grothendieck topology $\tau$ which is at 
least as fine as the \'etale topology.   

For $U\in \bar\site$, and $g\in \hat G$ we let
$gU=U\times_{\bar \F_q,g^{-1}}\bar \F_q$, so that for every sheaf 
$\f$ on $\bar\site$ we have $g^*\f(U)=\f(gU)$ and $g_*\f(U)=\f(g^{-1}U)$. 
We say that $\hat G$ acts of $\f$, if for every $g\in \hat G$ there is an 
isomorphism $\sigma_g:\f\to g^*\f$ satisfying
$\sigma_{gh}=h^*\sigma_g\circ \sigma_h$. For $f\in \f(U)$, we will
abbreviate $\sigma_g(f)\in \f(gU)$ by $gf$. 

Let $\hat S(U)\subseteq \hat G$ be the Galois group of the smallest
field extension $\F_{q^r}$ of $\F_q$ over which $U$ has a model $U'$, 
i.e. $U=U'\times_{\F_{q^r}}\bar\F_q$, and let $S(U)=\hat S(U)\cap G$.
If $\hat G$ acts on $\f$, then $\hat S(U)$ acts on $\f(U)$. In particular,
$\hat G$ acts on $\f(V\times_{\F_q}\bar\F_q)$ for every $V\in\site$.
We say that $\hat G$ acts continuously on $\f$, if for each
\'etale $U\in \hat\site$, $\hat S(U)$ acts
continuously on $\f(U)$ equipped with the discrete topology, 
i.e. if the map $\hat S(U)\times \f(U)\to \f(U)$ is continuous.
Let $\T_{\hat G}$ be the topos of sheaves on
$\bar\site$ equipped with a continuous action of $\hat G$.

\begin{lemma}\label{deligne}
a) If $\f$ is a sheaf on $\bar\site$, then $\hat G$ 
acts continuously on $\f$ if and only if 
$\colim_{H\subseteq \hat S(U)}\f(U)^H@>\sim >>\f(U)$ 
for every $U\in \bar\site$. The maps in the direct system 
are the natural inclusion maps.

b) There is an equivalence of categories between the category of sheaves 
on $\site$ and the category $\T_{\hat G}$.
\end{lemma}

\proof 
a) This is well-known.

b) This is Deligne [SGA 7 XIII, 1.1.3]. Explicitly, 
if $\pi:\Spec \bar \F_q\to \Spec \F_q$ is the structure map, then
the sheaf $\g$ on $\site$ corresponds to the sheaf $\pi^*\g$ on $\bar\site$, 
sending $U\to \bar X$ with model $U'$ over $\F_{q^r}$ to
$\colim_m \g(U'\times_{\F_{q^r}}\F_{q^{rm}})$.
The actions of $\hat G/m\hat G$ on
$\g(U'\times_{\F_{q^r}}\F_{q^{rm}})$ are compatible and give
an action of $\hat G$ on the colimit. Conversely, a sheaf $\f$
in $\T_{\hat G}$ corresponds to the sheaf $\pi_*^{\hat G}\f$ on $\site$, 
sending $V$ to $\f(V\times_{\F_q}\bar \F_q)^{\hat G}$.
\proofend

In \cite{licht}, Lichtenbaum defines the Weil-\'etale topology on
the small \'etale site of a scheme $X$ of finite type over $\F_q$.
He shows that a Weil-\'etale sheaf is equivalent to an
\'etale sheaf on $\bar X$ together with a $G$-action, where $n\in
G$ acts on $\bar X$ via $\varphi^n$ and on $\f$ via
$\sigma_\varphi^n:\f(U)\to \f(\varphi^nU)$.
In accordance with Lichtenbaum's definition, we let
$\T_G$ be the topos of sheaves on $\bar\site$ equipped with an
action of $G$, and call it the Weil-\'etale topos.

\begin{lemma}\label{adjoints}
The forgetful functor from abelian groups of $\T_G$ to sheaves of
abelian groups on $\bar\site$ has an exact left adjoint and a right adjoint.
In particular, it preserves injectives, and $\T_G$ has enough injectives.
\end{lemma}

\proof The left adjoint is $\f\to \oplus_{g\in G}g^*\f$ and the
right adjoint is $\f\to \prod_{g\in G}g^*\f$. In both cases, the
action of $G$ is the shift functor. A map
$\alpha\in\Hom_{\bar\site}(\f,\g)$, corresponds to the $G$-invariant
map $\Hom_{\T_G}(\oplus_{g\in G}g^*\f,\g)$ which on the 
summand indexed by $g$ is the composition $g^*\f@>g^*(\alpha)>>
g^*\g@>\sigma_g^{-1} >> \g$. The right adjoint case is analog.
Since $g^*$ and coproducts are exact, the left adjoint is exact 
and hence the forgetful functor preserves injectives. On the other
hand, given a sheaf $\f$ in
$\T_G$, we can embed it into an injective \'etale sheaf $\cal I$
on $\bar\site$. This gives rise to a $G$-invariant injection of $\f$
into the sheaf
$\prod_{g\in G}g^*\cal I$, which is injective in $\T_G$. 
\proofend

Recall that a morphism of topoi $\alpha: \T\to {\cal S}$ 
is a pair of adjoint functors $\alpha^*\vdash \alpha_*$ such
that $\alpha^*$ commutes with finite limits.

\begin{prop}
There is a morphism of topoi
$\gamma:\T_{G}\to \T_{\hat G}$. The functor
$\gamma^*$ is the forgetful functor, and
$$\gamma_*\f(U)=\colim_{H\subseteq S(U)}\f(U)^{H},$$
where $H$ runs through the subgroups of finite index in $G$ which
are contained in $S(U)$. In particular, $\gamma_*$
is left exact and preserves injectives. The adjoint transformation
$\id \to \gamma_*\gamma^*$ is an isomorphism.
\end{prop}

\proof Since the invariant functor is left exact, $\gamma_*$ of a
sheaf is sheaf. The action of $\hat G$ on
$\gamma_*\f$ is given as follows. Given
$U$ and $H\subseteq S(U)$,
$g\in \hat G$ acts as $\sigma_g=\sigma_\varphi^i: \f(U)^H \to
\f(gU)^H$, if $g\equiv \varphi^i\mod H$. It is easy to check that
this is compatible with the inclusion $\f(U)^H\hookrightarrow
\f(U)^{H'}$ for $H'\subseteq H\subseteq S(U)$, and hence induces
an action of $\hat G$ on the colimit.

Let $\f$ be a sheaf with $G$-action and $\g$ be a sheaf with
continuous $\hat G$-action. Then $\g\cong\colim_{H}\g^H$, and the
map
\begin{align*}
\Hom_G(\gamma^*\g,\f)&\to  \Hom_{\hat G}(\g,\gamma_*\f)\\
\alpha &\to \colim_H \alpha|_{\g^H}
\end{align*}
is an isomorphism with inverse "composition with the adjoint
inclusion $\gamma^*\gamma_*\f \to \f$". The fact
$\f@>\sim>> \gamma_*\gamma^*\f$ follows from the explicit
description of $\gamma_*$ and $\gamma^*$.
\proofend

Since subgroups $H\subseteq \hat S(U)$ are cofinal in
the set of all subgroups of finite index of $\hat G$, we will
write by abuse of notation
$\gamma_* \f =\colim_{H} \f^H$, remembering that even though
not every $\f^H$ is defined, the colimit is.

\section{The functor $R\gamma_*$}
Given two sheaves $\f$ and $\g$ in $\T_G$, the sheaf
$\HOM(\g,\f)$ is equipped with a $G$-action by $f^g= \sigma_g\circ
f \circ \sigma_g^{-1}$. Then
$$\gamma_*\HOM(\g,\f)(U)=
\colim_{H\subseteq S(U)} \Hom(\g|_{U},\f|_{U})^H,$$ where the
latter are the homomorphisms which are compatible with the action
of $H$. If $\g=A$ is constant, then by adjointness of global
section and constant sheaf functor we have 
$\Hom_{Shv}(A,\f|_{U})^H\cong \Hom_{\Ab}(A,\f(U))^H$, 
and the formula simplifies to
$\gamma_*\HOM(A,\f)(U)= \colim_{H\subseteq S(U)} \Hom(A,\f(U))^H$.

If $mG\subseteq S(U)$ and $i\in \Z/m$, then $\varphi^{\bar i}U$
does not depend on the representative $\bar i\in \Z$ of $i$, and we 
simply write $\varphi^iU$. We denote the $i$th summand of 
$f\in \oplus_{i\in \Z/m} \f(\varphi^iU)$ by $f^{(i)}$

\begin{lemma}
\label{coind} Let $\f$ in $\T_G$, $U\in \bar\site$, and
$H=mG\subseteq S(U) \subseteq G$.

a) If $H$ acts on $\Hom(\Z[G],\f(U))\cong\HOM(\Z[G],\f)(U)$
via
$F\mapsto h\circ F\circ h^{-1}$, then there are isomorphisms
\begin{align*}
\Z[G] \otimes_H \f(U)&@>\alpha >\sim > \Hom(\Z[G],\f(U))^H ,\\
\alpha(s\otimes f)(g)&=
\begin{cases}
gsf &{\text if\;} gs\in H;\\
0 &\text{otherwise.}
\end{cases}\\
\Z[G] \otimes_H \f(U)&@>\beta>\sim> \bigoplus_{i\in\Z/m}\f(\varphi^iU),\\
\beta(\varphi^a \otimes f)^{(i)}&=
\begin{cases}
\varphi^a f &{\text if\;} i\equiv a\mod m;\\
0 &\text{otherwise.}
\end{cases}
\end{align*}

b) The action of $\varphi\in G$ on $\Hom(\Z[G],\f(U))^H$ via 
multiplication on $\Z[G]$ corresponds under $\beta\alpha^{-1}$ to
the automorphism $\zeta(f)^{(i)}=\varphi f^{(i-1)}$ of 
$\bigoplus_{i\in \Z/m}\f(\varphi^iU)$.

c) If $g\in \hat G/m\hat G$ satisfies $g\equiv\varphi^a \mod m\hat G$, 
then the map $\Hom(\Z[G],\f(U))^H\to \Hom(\Z[G],\f(gU))^H$,
$F\mapsto g\circ F\circ g^{-1}$
corresponds under $\beta\alpha^{-1}$ to the cyclic permutation
$$\tau_g: \bigoplus_{i\in \Z/m}\f(\varphi^iU)\to \bigoplus_{i\in
\Z/m}\f(\varphi^{i}gU), \quad 
\tau_g(f)^{(i)}=f^{(i+a)}.$$

d) Given a second subgroup $H'=mnG \subseteq H\subseteq S(U)$, the
inclusion of fixed points $\Hom(\Z[G],\f(U))^H\hookrightarrow
\Hom(\Z[G],\f(U))^{H'} $ corresponds under $\beta\alpha^{-1}$ to the map
$$ \delta_m^n:\bigoplus_{i\in \Z/m}\f(\varphi^iU) \to \bigoplus_{j\in
\Z/mn}\f(\varphi^jU),\quad 
\delta^n_m(f)^{(j)}= f^{(j\mod m)}.$$

e) The map $\delta_m^n$ is compatible with the action of $G$ and
of $\hat G$ given in b) and c), respectively.
\end{lemma}

\proof a) This is an easy verification.

b) Consider the action of $G$ on $\Z[G] \otimes_H\f(U)$ by left
multiplication on $\Z[G]$. Then it is easy to verify that the
three actions are compatible with $\alpha$ and $\beta$.

c) The conjugation map is well defined, because since $F$ is
$H$-invariant, we have for $h\in H$,
$g\circ h\circ F\circ h^{-1}\circ g^{-1}=g\circ F\circ g^{-1}$.
If $g\equiv\varphi^a\mod m\hat G$, then
$gU=\varphi^aU$, so that $\f(\varphi^igU)=\f(\varphi^{i+a}U)$. 
It is easy to check that the conjugation map
$F\mapsto g\circ F\circ g^{-1}$ corresponds under $\alpha$
to the map
\begin{align*}
\Z[G] \otimes_H\f(U)&\to \Z[G] \otimes_H\f(gU)\\
s\otimes f&\mapsto g^{-1}s\otimes gf
\end{align*}
and this corresponds to $\tau$ under $\beta$.

d) Let $u$ be the inclusion of $H$-invariant maps into
$H'$-invariant maps, and let $v:\Z[G] \otimes_H \f(U) \to \Z[G]
\otimes_{H'} \f(U) $ be the map $s\otimes f \mapsto
\sum_{j=0}^{n-1} s\varphi^{jm}\otimes \varphi^{-jm}f$. It is easy
to check that $u\circ\alpha=\alpha\circ v$ and $\delta_m^n\circ
\beta=\beta\circ v$.

e) This is clear because the actions of $G$ and $\hat G$ are
compatible with the inclusion map, and with $\alpha$ and
$\beta$. Explicitly,
$$\delta_m^n\zeta(f)^{(j)}=\zeta(f)^{(j\mod m)}=
\varphi f^{(j-1 \mod m)}=\varphi \delta_m^n(f)^{(j-1)}
=\zeta\delta_m^n(f)^{(j)},$$
and for $g\in \hat G$ with $g\equiv \varphi^a\mod mn\hat G$,
$$ \delta^n_m\tau_g(f)^{(j)}=\tau_g(f)^{(j\mod m)} =f^{(j+a\mod m)}
=\delta_m^n(f)^{(j+a)}=\tau_g\delta_m^n(f)^{(i)}.$$
\proofend

Consider the presheaf
$$\fun(\f):U\mapsto \colim_{mG\subseteq S(U)}
\bigoplus_{i\in \Z/m}\f(\varphi^iU),$$ 
where the index set is ordered by divisibility, and the maps in the direct 
system are the maps $\delta_m^n$. The presheaf $\fun(\f)$ is a 
sheaf, because filtered direct limits and direct sums are left exact.
Moreover, the action of $g\in \hat G$ is compatible with $\delta_m^n$, 
so that we get an action of $\hat G$ on $\fun(\f)$.

\begin{lemma}
\label{vanish} The functor $\f\mapsto \gamma_* \HOM(\Z[G],-)$ from
$\T_G$ to $\T_{\hat G}$ is exact. In particular, the derived
functors $R^i\colim_H \HOM(\Z[G],- )^H$ are zero for $i>0$.
\end{lemma}

\proof By Lemma \ref{coind}, $\colim_m \HOM(\Z[G],\f)^{mG} \cong
\colim_m \bigoplus_{i\in \Z/m}(\varphi^{i})^*\f$. Now the functor
$\f\mapsto (\varphi^{i})^*\f$ and filtered colimits of sheaves
are exact. \proofend

\begin{theorem}
\label{resolution} Let $\f$ be a sheaf in $\T_G$. Then the complex
$R\gamma_*\f$ is quasi-isomorphic to the complex of sheaves of
continuous $\hat G$-modules sending $U\in \bar\site$ to
\begin{equation}\label{twoc}
\fun(\f)(U)@>t-1>>\fun(\f)(U).
\end{equation}
Here $(tf)^{(i)}=\varphi f^{(i-1)}$ and $g\in \hat G$ with
$\varphi^a\equiv g\mod m\hat G$ acts as $(gf)^{(i)}=f^{(i+a)}$
on $\oplus_{i\in \Z/m}\f(\varphi^iU)$.
\end{theorem}

\proof 
Let $P_\cdot$ be the free resolution
$0 \to \Z[G] @>t-1>> \Z[G]  \to 0$ of the constant sheaf $\Z$,
and let $ \f \to I^\cdot $ be an injective resolution. Then
$R\gamma_*\f$ is quasi-isomorphic to
$\colim_m(I^\cdot)^{mG}\cong \colim_m \HOM(\Z,I^\cdot)^{mG}
\cong \colim_m \HOM(P_\cdot,I^\cdot)^{mG}$. If we take vertical
cohomology in the latter double complex, we get complexes
$$R^b\big( \colim_m\HOM(\Z[G],-)^{mG} \big)(\f) @>t-1 >>
R^b\big( \colim_m\HOM(\Z[G],-)^{mG} \big)(\f)$$ concentrated in
degree $a=0,1$ for each $b$. But by Lemma \ref{vanish} the derived
functors vanish for $b>0$, and the double complex is
quasi-isomorphic to
$$ \colim_m \HOM(\Z[G],\f)^{mG} @>t-1>>
\colim_m \HOM(\Z[G],\f)^{mG} ,$$ where the map $t$ is given in
Lemma \ref{coind} b) and the $\hat G$-action in Lemma \ref{coind} c). 
By Lemma \ref{coind} a), this complex is isomorphic to the complex of the
theorem. 
\proofend

\section{$\gamma_*\f$, $R^1\gamma_*\f$, and $-\cup e$}

To calculate the cohomology sheaves of $R\gamma_*\f$ explicitly,
let
$N_m^n:\f(U)\to \f(U)$ be the map $f\mapsto
\sum_{l=0}^{n-1}\varphi^{-lm}f$ for
$mG\subseteq S(U)$. This descends to a map $N_m^n:\f(U)_{m\Z}\to
\f(U)_{mn\Z}$, because
$$N_m^n(\varphi^mf)=\sum_{l=0}^{n-1}\varphi^{-(l-1)m}f =
N_m^n(f)+\varphi^m f-\varphi^{-(n-1)m}f
=N_m^n(f)+(1-\varphi^{-nm})\varphi^{m}f.$$

\begin{prop}
\label{maps}
a) Let $\f$ be a sheaf in $\T_{G}$. Then there is an exact 
sequence
$$0\to \colim_m \f^{mG} @>\Delta>> \fun(\f) @>t-1 >>\fun(\f)
@>S>> \colim_{m,N_m^n}\f_{mG}\to 0.$$

b) If $\f$ is a sheaf of $\Q$-vector spaces, then $\Delta$ is split, hence
$$R\gamma_*\f\cong \gamma_*\f\oplus R^1\gamma_*\f[-1].$$
\end{prop}

\proof a) We have to calculate $\gamma_*\f=\ker t-1$
and $R^1\gamma_*\f= \coker t-1$.
To construct the isomorphism $\Delta:\colim_m \f^{mG}\to \ker t-1$,
let $\f$ in $\T_{G}$, fix
$U\in \bar\site$, and let $mG\subseteq S(U)$. Define
$$\Delta_m : \f(U)^{mG} \to \bigoplus_{i\in \Z/m}\f(\varphi^iU),
\quad  \Delta_m(c)^{(i)}=\varphi^{\bar i}c ,$$
where $\bar i\in \Z$ is the representative of $i\in \Z/m$
with $0\leq \bar i\leq m-1$. 
The image of  $\Delta_m$ is contained in
the kernel of $t-1$, because
$$((t-1)\Delta_m(c))^{(i)} =
\varphi\Delta_m(c)^{(i-1)}-\Delta_m(c)^{(i)}=
\varphi^{\overline{i-1}-\bar i+1}f=0.$$
Clearly $\Delta_m$ is injective, and if
$f\in \bigoplus_{i\in \Z/m}\f(\varphi^iU)$ satisfies
$(t-1)f =0$, then $f^{(i)}=\varphi^{\bar i}f^{(0)}$,
so that $\Delta_m(f^{(0)})=f$, and
$\Delta_m$ is an isomorphism to the kernel of $t-1$. It is easy
to check that $\delta^n_m\circ\Delta_m = \Delta_{mn}$, hence we get a map
\begin{equation}\label{delta}
\Delta: \colim_m  \f(U)^{mG} \to \fun(\f)(U)
\end{equation}
which is an isomorphism to the kernel of $t-1$.

To construct the isomorphism $S:\coker t-1 \to \colim_{m,N_m^n}\f_{mG}$, 
consider the map
$$S_m : \bigoplus_{i\in \Z/m}\f(\varphi^iU) \to \f(U), \quad 
f\mapsto \sum_{i\in \Z/m}\varphi^{-\bar i}f^{(i)}.$$
We have
\begin{multline*}
S_{mn}(\delta^n_m (f)) =
\sum_{j\in \Z/mn}\varphi^{-\bar j}\delta_m^n(f)^{(j)}
=\sum_{j\in \Z/mn}\varphi^{-\bar j}f^{(j\mod m)}\\
=\sum_{l=0}^{n-1}\sum_{i\in \Z/m}\varphi^{-\bar i-lm}f^{(i)}
=N_m^n(\sum_{i\in \Z/m}\varphi^{-\bar i}f^{(i)})=N_m^n(S_m(f)).
\end{multline*}
hence a surjective map
$$\tilde S:\fun(\f)(U)
\to \colim_{m,N^n_m} \f(U).$$ 
Since
\begin{multline*} S_m((t-1)f)=
\sum_{i\in \Z/m}\varphi^{\bar i}(\varphi f^{(i-1)}-f^{(i)})\\
=\varphi f^{(m-1)}-\varphi^{-m+1}f^{(m-1)}=
(1-\varphi^{-m})(\varphi f^{(m-1)}),
\end{multline*}
this map induces a map
\begin{equation}\label{sss}
S:\Big( \colim_m  \bigoplus_{i\in \Z/m}\f(\varphi^iU)\Big)/t-1
\to \colim_{m,N^n_m} \f(U)_{mG},
\end{equation}
and we claim that $S$ is an isomorphism. Indeed, for 
$c\in \f(\varphi^jU)$ define 
$R_j(c)\in\bigoplus_{i\in \Z/m}\f(\varphi^iU)$ by
$$ R_i(c)^{(l)}=
\begin{cases}
\varphi^{l-\bar j}c &l<\bar j;\\
0& l\geq \bar j.
\end{cases}$$
Then 
$$((t-1)R_i(c))^{(i)}=
\begin{cases}
c& i=j;\\ 
-\varphi^{-\bar j}c &i=0;\\
0&\text{otherwise}.
\end{cases}$$
Hence 
$$(f-(t-1)\sum_{i\in \Z/m} R_i(f^{(i)}))^{(l)}
=\begin{cases}
\sum_{i\in \Z/m}\varphi^{-\bar i}f^{(i)}&l=0\\
0 &l\not=0.
\end{cases}$$
We conclude that $S_m(f) =\sum_i \varphi^{-\bar i}f^{(i)}=0$
implies that $f$ is in the image of $t-1$.

b) This follows because the map
$S'=\colim_m\frac{1}{m}S_m:\fun(\f)(U)
\to \gamma_*\f(U)$ satisfies $S'\circ \Delta=\id$, hence we get a
quasi-isomorphism
$$\begin{CD}
\fun(\f)(U) @>t-1>> \fun(\f)(U)\\
@VS' VV @VSVV\\
\gamma_*\f(U)@>0>> R^1\gamma_*\f(U).
\end{CD}$$
\proofend

\begin{cor}
We have $\gamma_*\Z\cong \Z$, $R^1\gamma_*\Z\cong \Q$, and
$R\gamma_*\Q\cong \Q\oplus \Q[-1]$.
\end{cor}

\proof
The map which is multiplication by $\frac{1}{m}$ on the
copy of $\Z$ indexed by $m$ induces an isomorphism
$\colim_{m,N^n_m}\Z@>\sim>> \Q$.
\proofend

\medskip\noindent{\bf Examples.}
1) If a generator of $G$ acts on the constant sheaf $\f= \Q$ as
multiplication by $r\not=\pm 1$, then $R\gamma_*\f=0$. More
generally, for any constant sheaf $\f$ of rational vector spaces,
$\gamma_*\f$ is the largest subspace on which $\varphi$ acts as
multiplication by some root of unity, and $R^1\gamma_*\f$ is the
largest quotient space on which $\varphi$ acts as multiplication
by some root of unity.

\smallskip

2) Let $\f$ be the sheaf $\oplus_{i\in G} A$ for an abelian group
$A$, where $G$ acts by shifting the factors. Then $\gamma_*\f=0$,
whereas $R^1\gamma_*\f\subseteq \prod_i A$ consists of elements
which are periodic for some period. Indeed, for an element of
$\f^{mG}$, the entries in the sum are $m$-periodic, hence they
must be zero. On the other hand, $\f_{mG}\cong \oplus_{i\in
\Z/m}A$ and under this isomorphism the map $\f_{mG}\to \f_{mnG}$
is the $n$-fold concatenation map.

\medskip

Consider the extension
$e\in \Ext^1_G(\Z,\Z)$ of $G$-modules, given by $N=\Z\oplus \Z$ as
an abelian group, and $g\cdot(r,s)=(r+as,s)$ for $g\in G$.

\begin{lemma}
The extension $e$ is a generator of
$\Ext^1_G(\Z,\Z)\cong \Z$.
\end{lemma}

\proof First note that $\Ext^1_G(\Z,\Z)\cong H^1(G,\Z)\cong \Z$.
Any extension $E\in \Ext^1_G(\Z,\Z)$, is isomorphic to
$\Z\oplus \Z$ as an abelian group. The action of a generator of
$G$ is given by a matrix of the form
${1\,n\choose 0\,1}$ for $n$ an integer. 
The formula for addition of extensions classes
then shows that $E=n\cdot e\in \Ext^1_G(\Z,\Z)$. 
\proofend

For every complex of Weil-\'etale sheaves $\f^\cdot$, the
connecting homomorphism of the distinguished triangle
$$ \f^\cdot \to \f^\cdot\otimes N \to \f^\cdot @>\beta>>  \f^\cdot[1]$$
induces cup product with $e$ on cohomology (up to sign).
Indeed, $e$ is the image of $\id$ under the induced map
$\Hom_G(\Z,\Z)@>\beta >> \Ext^1_G(\Z,\Z)$, and cup 
product is compatible with the Bockstein homomorphism, so that
$e\cup x=\beta(\id)\cup x=\beta(\id\cup x)=\beta(x)$.

\begin{prop}
\label{cupe} Let $\f$ be a Weil-\'etale sheaf. Then under the
identification of \eqref{twoc}, the cup product
with $e$ on $R\gamma_*\f$ is induced by the following map 
$R\gamma_*\f\to R\gamma_*\f[1]$ in the derived category of \'etale sheaves
$$\begin{CD}
R\gamma_*\f @EEE @EEE \fun(\f)@>t-1>> \fun(\f) \\
@VeVV  @III @V\id VV  \\
R\gamma_*\f[1] @EEE \fun(\f)@>t-1>> \fun(\f)
\end{CD}$$
\end{prop}

\proof Cup product with $e$ is induced by the following vertical 
map of double complexes
$$\begin{CD}
R\gamma_*\f @>>> R\gamma_*(\f\otimes N)\\
@V\id VV \\
R\gamma_*\f.
\end{CD}$$
Let $\alpha:R\gamma_*\f\to R\gamma_*\f[1]$ be the map in the statement
of the proposition, except that we replace the vertical identity
map $\fun(\f)\to \fun(\f)$ by the cyclic permutation map 
$\fun(\f)@>t>> \fun(\f)$. In the
derived category, $\alpha$ is quasi-isomorphic to the map
$$\begin{CD}
R\gamma_*\f @>>> cone(\alpha)\\
@V\id VV  \\
R\gamma_*\f.
\end{CD}$$
Cup product with $e$ and $\alpha$ agree in the derived category,
because in both diagrams, the upper row is the following double complex
$$\begin{CD}
\fun(\f)@>i_1 >> \fun(\f)\oplus \fun(\f) \\
@V t-1VV @V {t-1\,\,\,\,\,t\choose 0\,\,\,\,\,t-1}VV\\
\fun(\f)@>i_1>>\fun(\f)\oplus \fun(\f).
\end{CD}$$
Finally, the map $\alpha$ and the map of the proposition
are homotopic via the chain homotopy $h:R\gamma_*\f\to
R\gamma_*\f$, which is the identity map in degree $0$ and the zero
map in degree $1$. 
\proofend

As a consequence, we see that under 
the identification of Proposition \ref{maps} a), the cup product with 
$e$ induces the colimit of the canonical maps 
$\f^{mG}\to \f_{mG}$, $f\mapsto mf$ on cohomology sheaves.
Indeed, cup product with $e$ is the composition
$$ \gamma_*\f\cong \colim_m \f^{mG}@>\Delta >> \fun(\f) 
@>S>> \colim_{m,N_m^n}\f_{mG}\cong R^1\gamma_*\f,$$ 
and $S_m\circ \Delta_m(f)=mf$.

\section{The functor $R\gamma_*\gamma^*$}

\begin{theorem}\label{maincup} 
Let $\g^\cdot $ be a complex of sheaves in
$\T_{\hat G}$. Then there is a quasi-isomorphism
$$R\gamma_*\Z\otimes^L \g^\cdot \cong  R\gamma_*(\gamma^*\g^\cdot).$$
\end{theorem}

\proof The complex $R\gamma_*\Z$ of Theorem \ref{resolution} 
consists of {\it flat} sheaves. Hence we get a
quasi-isomorphism of complexes
$$ R\gamma_*\Z\otimes^L\g^\cdot (U) \cong
\colim_m \bigoplus_{i\in \Z/m}\g^\cdot (U)  @>t-1>> \colim_m
\bigoplus_{i\in \Z/m}\g^\cdot (U).$$ Here $t$ is the cyclic
permutation of the summands.
On the other hand, we have a quasi-isomorphism
$$R\gamma_*(\gamma^*\g^\cdot ) (U) \cong
\colim_m  \bigoplus_{i\in \Z/m}\g^\cdot(\varphi^iU) @>t-1>>
\colim_m \bigoplus_{i\in \Z/m}\g^\cdot (\varphi^iU).
$$
Here $t$ acts as in Theorem \ref{resolution}. We claim that the
following maps induce a quasi-isomorphism of the two complexes (we
only write one of the two identical terms of the complexes):
\begin{equation}\label{yuio}
\colim_m \bigoplus_{i\in \Z/m}\g^\cdot (U) \supseteq
\colim_m \bigoplus_{i\in \Z/m}\g^\cdot (U)^{mG}\hookrightarrow
\colim_m \bigoplus_{i\in \Z/m}\g^\cdot(\varphi^iU) .
\end{equation}
The restriction of $\delta_m^n$ to the middle subgroup is defined
because $\g^j(U)^{mG}\subseteq \g^j(U)^{mnG}$.
The right map $\nu$ is given by $\nu(f)^{(i)}=\varphi^{\bar i}f^{(i)}$. 
It is easy to verify that $\nu$ is compatible with 
the maps $\delta_m^n$, and the action of $g\in \hat G$. It also
commutes with the action of $t$, because
$(t\nu(f))^{(0)}=\varphi\nu(f)^{(m-1)}=\varphi^mf^{(m-1)}=f^{(m-1)}
=(tf)^{(0)}=\nu(tf)^{(0)}$ (here we see why we need to restrict $\nu$
to the the intermediate group).

The two inclusions are in fact bijections.
Indeed, because $\g^j(U)$ is a continuous $\hat G$-module, every 
element $x\in \g^j(U)$ in the $m$th term of the
colimit on the left is contained in $\g^j(U)^{nG}$ for some
$n$. Then $\delta_m^n(x)$ will be in the image of the inclusion map in the
$mn$th term of the colimit. The same argument works for $\nu$. 
\proofend

\begin{cor}
\label{rational} If $\g^\cdot $ is a complex of sheaves in
$\T_{\hat G}$, then there is a distinguished triangle
\begin{equation}\label{oror}
\g^\cdot \to R\gamma_*\gamma^*\g^\cdot \to \g^\cdot\otimes \Q[-1]
@>\delta>> \g^\cdot[1].
\end{equation}
In particular, if $\g^\cdot$ is a complex with torsion cohomology sheaves, then
$$\g^\cdot\cong R\gamma_*\gamma^*\g^\cdot.$$
If $\g^\cdot$ is a complex with $\Q$-vector spaces as cohomology sheaves, then
$$R\gamma_*\gamma^*\g^\cdot\cong \g^\cdot\oplus\g^\cdot[-1],$$
and under this isomorphism, cup product with $e$ is given by multiplication 
by the matrix ${0\,0\choose 1\, 0}$.
\end{cor}

\proof 
Everything is clear from the previous section except the
last statement. It is easy to check that the quasi-isomorphism
\eqref{yuio} is compatible with the description of the cup product
with $e$ in Proposition \ref{cupe}, hence it suffices to 
calculate the action of cup product with $e$ on $R\gamma_*\Z$.
By the remark after Proposition \ref{cupe}, this is induced by the 
composition
$$\Z\cong \colim_m \Z @>\Delta>>  \fun(\Z) @>S >> 
\colim_{m,N_m^n}\Z\cong \Q,$$
which is the inclusion map.
\proofend

We now calculate the map $\delta$. Consider the extension
$\epsilon \in \Ext^1_{\hat G}(\Q/\Z ,\Q/\Z)$, which is $\bar
N=\Q/\Z\oplus \Q/\Z$ as an abelian group, and $g\in \hat G$
acts via $g\cdot(x,y)=(x+gy,y)$. The image of $\epsilon$
under the composition of the canonical projection and the 
Bockstein homomorphism
$$\Ext^1_{\hat G}(\Q/\Z ,\Q/\Z)
@>>> \Ext^1_{\hat G}(\Q,\Q/\Z)
@>\beta >> \Ext^2_{\hat G}(\Q,\Z)$$
is calculated by the following pull-back commutative diagram
$$\begin{CD}
0@>>> \Z @>>> \Q @>\xi >> M @>>> \Q @>>> 0\\
@III @VVV @VVV @| @| \\
@EEE 0 @>>> \Q/\Z  @>>> M @>>> \Q @>>> 0\\
@III @III @| @VVV @VVV \\
@EEE 0 @>>> \Q/\Z  @>>> \bar N  @>>> \Q/\Z @>>> 0.
\end{CD}$$
Note that the cup product
with $\epsilon$ is the boundary map induced by the lower 
sequence. These extensions have been studied by Kahn
\cite[Def. 4.1]{kahn}. He denotes the complex
$\Q@>\xi >> M$ by $\Z^c$, where $\Q$ is in degree $0$, $M$
is in degree $1$. 

\begin{theorem}
\label{kahnequal}
There is a quasi-isomorphism of complexes of $\hat G$-modules
$R\gamma_* \Z \cong \Z^c$.
In particular, the boundary map $\delta$ in \eqref{oror} is the composition
\begin{equation}
\label{kahn} \Q[-1] \to \Q/\Z[-1]\ @>\epsilon>> \Q/\Z[0] @>\beta >>
\Z[1] .
\end{equation}
\end{theorem}

\proof Consider the map
\begin{align*} 
\mu_m:\bigoplus_{i\in \Z/m}\Z &\to M=\Q/\Z\oplus \Q,\\
\textstyle x&\mapsto \sum_{i\in \Z/m} \textstyle
(\frac{1}{2}x^{(i)}+\frac{ i}{m}x^{(i)}+\Z,\frac{-1}{m}x^{(i)}).
\end{align*}
This is a map of $\hat G$-modules, because for
$a \in \hat G/m\hat G$, we have
\begin{multline*}\textstyle \mu_m(ax)=
\sum_{i\in \Z/m}\textstyle (\frac{1}{2}x^{(i+a)}+\frac{i}{m}x^{(i+a)}+\Z,
\frac{-1}{m}x^{(i+a)})\\
=\sum_{i\in \Z/m}\textstyle (\frac{1}{2}x^{(i)}+\frac{i-a}{m}x^{(i)}+\Z,
\frac{-1}{m}x^{(i)})
=a\mu_m(x).
\end{multline*}
The maps $\mu_m$ are compatible with $\delta_m^n$ (here we need
the correcting summand $\frac{x^{(i)}}{2}$, and use the identity
$\sum_{l=0}^{n-1}\frac{i+lm}{nm}=\frac{i}{m}+\frac{n-1}{2}$):
\begin{multline*}\textstyle
\mu_{mn}(\delta^n_m(x))=
\sum_{j\in \Z/mn} (\frac{1}{2}x^{(j\mod m)}+\frac{l}{mn}x^{(j\mod m)}+\Z,
\frac{-1}{mn}x^{(j\mod m)})\\
\textstyle
=\sum_{i\in \Z/m} (\frac{n}{2}x^{(i)}+
\sum_{l=0}^{n-1}\frac{i+lm}{mn}x^{(i)}+\Z,
\frac{-n}{mn}x^{(i)})
=\mu_m(x).
\end{multline*}
Hence we get a map $\mu: \colim_m\oplus_{i\in \Z/m}\Z \to M$,
which is compatible with the action of $\hat G$. Consider the
following diagram of maps of $\hat G$-modules with exact rows,
where $S'=\colim_m\frac{1}{m}S_m$.
$$\begin{CD}
\Z @>\Delta >> \colim_m \bigoplus_{i\in \Z/m}\Z @>t-1
>> \colim_m \bigoplus_{i\in\Z/m} \Z @>S'>> \Q \\
@| @VS' VV @V \mu VV @| \\
\Z @>>> \Q @>\xi>> M @>p>> \Q
\end{CD}$$
It is easy to see that the outer squares commute. On the other hand,
if $x\in \oplus_{i\in \Z/m}\Z$, then
\begin{multline*}\textstyle
\mu_m(t-1)x=\sum_{i\in \Z/m} 
(\frac{1}{2}(x^{(i-1)}-x^{(i)})+\frac{i}{m}(x^{(i-1)}-x^{(i)})+\Z,
\frac{-1}{m}(x^{(i-1)}-x^{(i)}))\\
\textstyle
=\sum_{i\in \Z/m}(\frac{1}{m}x^{(i)}+\Z,0)
=\xi S'(x).
\end{multline*}
\proofend

\section{Weil-\'etale cohomology}
Let $\Ab$ be the category of abelian groups, and $\Mod_G$ and
$\Mod_{\hat G}$ are the categories of $G$-modules and continuous
$\hat G$-modules (for the discrete topology), respectively.
Consider the following commutative diagram of functors, where
$\Gamma_{\bar X}(\f)=\f(\bar X)$,
$$\begin{CD}
\T_G @>\gamma_*>>  \T_{\hat G} \\
@V \Gamma_{\bar X} VV @V \Gamma_{\bar X} VV \\
\Mod_G  @> \gamma_* >> \Mod_{\hat G}  \\
@V \Gamma_G VV @V \Gamma_{\hat G} VV \\
\Ab @= \Ab
\end{CD}$$
For a complex of sheaves $\f^\cdot$ in $\T_G$, we define the
derived functors
\begin{align*}
H^i_W(X,\f^\cdot)&=R^i(\Gamma_G\circ\Gamma_{\bar X})(\f^\cdot)\in \Ab\\
H^i_{\hat G}(\bar X,\f^\cdot)&=
R^i(\gamma_*\circ \Gamma_{\bar X})(\f^\cdot)\in \Mod_{\hat G}
\end{align*}
Following Lichtenbaum, we call the first groups Weil-\'etale
cohomology groups. 

\begin{lemma}\label{uiop}
For a complex of Weil-\'etale sheaves $\f^\cdot$,
the underlying abelian group of the derived functors 
$R^i\Gamma_{\bar X}(\f^\cdot)\in \Mod_G$ agree with the
usual \'etale cohomology groups $H^i_\et(\bar X,\f^\cdot)$.
\end{lemma}

\proof
The horizontal forgetful functors in the diagram
$$\begin{CD}
\T_G @>F>>  \bar\site \\
@V \Gamma_{\bar X} VV @V \Gamma_{\bar X} VV \\
\Mod_G  @> F >> \Ab  \\
\end{CD}$$
are exact, and all functors
involved have an exact left adjoint by Lemma \ref{adjoints}, hence
preserve injectives. Thus $F\circ R\Gamma_{\bar X}=R\Gamma_{\bar X}\circ F$.
\proofend

We get the following spectral sequences for
composition of functors
\begin{equation}
\begin{CD}\label{ss1} 
E_2^{s,t}=H^s_{cont}(\hat G,H^t_{\hat G}(\bar X,\f^\cdot))
@E\Rightarrow EE H^{s+t}_W(X,\f^\cdot)\\
@VVV @| \\
E_2^{s,t}=H^s(G,H^t_\et(\bar X,\f^\cdot)) @E\Rightarrow EE
H^{s+t}_W(X,\f^\cdot).
\end{CD}
\end{equation}
Since $G$ has cohomological dimension $1$, the latter spectral
sequence breaks up into short exact sequence
\begin{equation}
\label{descent} 0\to H^{t-1}_\et(\bar X,\f^\cdot)_G\to
H^t_W(X,\f^\cdot) \to H^t_\et(\bar X,\f^\cdot)^G \to 0 .
\end{equation}
The identity
$\Gamma_{\hat G}\circ\Gamma_{\bar X}\circ \gamma_*\cong 
\Gamma_G\circ \Gamma_{\bar X}$
gives another spectral sequence
$$ H^s_\et(X,R^t\gamma_*\f^\cdot)\Rightarrow H^{s+t}_W(X,\f^\cdot).$$

\begin{lemma}\label{multss}
a) For a $G$-module $A$, the cup product with $e\in H^1(G,\Z)$ induces 
the canonical map
$H^0(G,A)\cong A^G\to A_G\cong H^1(G,A)$ on cohomology.

b) Cup product with $e$ agrees with the composition
$$ H^t_W(X,\f^\cdot)@>>> H^t_\et(\bar X,\f^\cdot)^G 
@>can >> H^t_\et(\bar X,\f^\cdot)_G @>>> H^{t+1}_W(X,\f^\cdot).$$
\end{lemma}

\proof
a) This is proved as Proposition \ref{cupe}, see \cite[Lemma 1.2]{rapzink}. 

b) Because the Hochschild-Serre spectral sequence is multiplicative,
and $H^2(G,M)=0$ for every $G$-module $M$, the cup product
with $e\in H^1_W(\F_q,\Z)\cong H^1(G,\Z)= E^{1,0}_2(\F_q)$ on 
the abutment of \eqref{ss1} is induced by the cup product 
$E_2^{1,0}(\F_q)\times E_2^{0,t}(X)\to E_2^{1,t}(X)$. Hence the result
follows from a).
\proofend

\subsection{Continuous Weil-\'etale cohomology}
For pro-\'etale sheaves, Jannsen \cite{jannsen} defined continuous
cohomology as the derived functors of $\lim\circ \Gamma_X$. We are
introducing the analog for the Weil-\'etale topology.

\begin{lemma}
\label{lim} Limits in the categories $\T_G$ and $\T_{\hat G}$
exist. More precisely, $\lim_{G}$ agrees with the limit in the
category of \'etale sheaves on $\bar X$, and
$\lim_{\hat G} \cong \gamma_*{\lim_G} \gamma^* $.
In particular, $R\lim_{\hat G}\cong R\gamma_*R{\lim_G} \gamma^* $.
\end{lemma}

\proof 
The first statement follows from Lemma \ref{adjoints}.
Since $\gamma_*$ has a left adjoint, it commutes with limits, and
we get
$\gamma_*\lim_G\gamma^*\cong \lim_{\hat G}\gamma_*\gamma^*\cong \lim_{\hat G}$.
\proofend

Let $\T_G^\N $ and $\T_{\hat G}^\N$ be the category of inverse
systems, indexed by the natural numbers, of sheaves in $\T_G$ and
$\T_{\hat G}$, respectively. There is a commutative diagram of
functors
$$ \begin{CD}
\T_G^\N @> \gamma_* >> \T_{\hat G}^\N \\
@V \lim_G VV @V \lim_{\hat G}  VV \\
\T_G @>\gamma_* >> \T_{\hat G}.
\end{CD}$$
By Lemma \ref{deligne} b), the functor $\lim_{\hat G}$ on
$\T^\N_{\hat G}$ corresponds to the functor $\lim$ on
$\site$ under the identification of Lemma \ref{deligne}. 
For $\f_\cdot\in \T_{G}^\N$ a pro-system of sheaves of
$G$-modules on $\bar X$, we define derived functors
\begin{align*}
H^i_{W}(X,(\f_\cdot))&=
R^i(\Gamma_{G}\circ\lim_{G}\circ\Gamma_{\bar X})(\f_\cdot)\in \Ab\\
H^i_{\hat G}(\bar X,(\f_\cdot))&=
R^i(\gamma_*\circ\lim_{G}\circ\Gamma_{\bar X})(\f_\cdot)\in
\Mod_{\hat G}.
\end{align*}
The second groups are continuous $\hat G$-modules for the {\it discrete} 
topology. It follows as in Lemma \ref{uiop} that
the underlying abelian group of the derived functors
$R^i(\lim_{G}\circ\Gamma_{\bar X})(\f_\cdot)\in \Mod_G$
agrees with the usual continuous cohomology groups of Jannsen of
$\bar X$. There is a map of spectral sequences
\begin{equation}
\begin{CD}\label{moress} 
E_2^{s,t}=H^s_{cont}(\hat G,H^t_{\hat G}(\bar X,(\f_\cdot)))
@E\Rightarrow EE
H^{s+t}_{W}(X,(\f_\cdot))\\
@VVV @| \\
E_2^{s,t}=H^s(G,H^t(\bar X,(\f_\cdot))) @E\Rightarrow EE
H^{s+t}_{W}(X,(\f_\cdot)).
\end{CD}
\end{equation}

\begin{lemma}
For an inverse system $\g_\cdot$ of torsion sheaves in $\T_{\hat
G}^\N$, the cohomology groups $H^i_{W}(X,(\gamma^*\g_\cdot))$
agree with the continuous cohomology groups
$H^i(X,(\g_\cdot))$ of Jannsen on $X$, and the groups
$H^i_{\hat G}(\bar X,(\gamma^*\g_\cdot))$ agree with the
continuous cohomology groups
$H^i(\bar X,(\g_\cdot))$ of Jannsen on $\bar X$.
In particular,
$$H^i_{cont}(X,\Z_l(n)):=H^i(X,(\Z/l^\cdot(n))) \cong
H^i_{W}(X,(\Z/l^\cdot(n))).$$
\end{lemma}

\proof Note first that $\Gamma_G\circ \lim_G\circ \Gamma_{\bar
X}=\lim_{\hat G}\circ \Gamma_{\hat G}\circ \Gamma_{\bar X}\circ
\gamma_*$. Now $R\gamma_*\gamma^*\g_\cdot \cong \g_\cdot$ for
torsion sheaves, and the derived functors of $\lim_{\hat G}\circ
\Gamma_{\hat G}\circ \Gamma_{\bar X}$ are the continuous
cohomology groups in the sense of Jannsen. The second statement
follows similarly. \proofend

The upper spectral sequence \eqref{moress} differs from Jannsen's
spectral sequence \cite[Cor. 3.4]{jannsen}, for finitely generated
cohomology groups
$$ H^s_{cont}(\hat G,H^t(\bar X,(\g_\cdot))) \Rightarrow
H^{s+t}(X,(\g_\cdot)),$$ even though for \'etale torsion sheaves
the $E_2$-terms and the abutment agree. This is because Jannsen
considers the coefficients with the limit topology, whereas we
work with the discrete topology.

\begin{lemma}
\label{stalks} a) If $(\f_\cdot)$ is an inverse system of sheaves
in $\T_{G}$, then the sheaf $R^s\lim\f_\cdot $ is the sheaf
associated to the presheaf which sends $U\in \bar\site$ to
$H^s_\et(U,(\f_\cdot))$.

b) If $(\g_\cdot)$ is an inverse system of sheaves in $\T_{\hat
G}$, then the sheaf $R^s\lim_{\hat G}\g_\cdot$ is the sheaf
associated to the presheaf
$U\mapsto H^s_{\hat G}(U,(\gamma^*\g_\cdot))$.
\end{lemma}

\proof a) follows with the same proof as b) by erasing all
$\gamma_*$ and $\gamma^*$.

b) Let $\g_\cdot \to I_\cdot^*$ be an injective resolution of
inverse systems. Then $R^s\lim_{\hat G}(\g_\cdot)$ is by
definition ${\cal H}^s(\gamma_*\lim \gamma^* I_\cdot^*) = a{\cal
H}^s(i\gamma_*\lim\gamma^*  I_\cdot^*)$, where $a$ is the
sheafification functor and $i$ the inclusion of sheaves into
presheaves. On the other hand, for every pro-system of sheaves
$(\lim\g_\cdot )(U)= \lim(\g_\cdot (U))$, for every sheaf
$(\gamma_*\g)(U)=\gamma_*(\g(U))$ and $\gamma^*\g(U)=\g(U)$, hence
$${\cal H}^s(i\gamma_*\lim \gamma^* I_\cdot^*)(U)=
H^s(\gamma_*\lim I_\cdot^*(U)) =: H^s_{\hat
G}(U,(\gamma^*\g_\cdot)).$$
\proofend

\section{Motivic cohomology}
From now on we will assume that $\site$ is the category
of smooth schemes of finite type over $\F_q$, or the small \'etale site
of a smooth scheme over $\F_q$ (Weil-\'etale cohomology does not
have good properties for non-smooth schemes).
For $n\geq 0$, let $\Z(n)$ be the (\'etale) motivic complex of 
Voevodsky \cite{susvoe}; for example $\Z(0)\cong \Z$
and $\Z(1)\cong {\mathbb G}_m[-1]$.
For an abelian group $A$, we define $A(n)$ to be $\Z(n)\otimes A$. 
In order to make our formulas work in general, we also define
$$ \Z(n)=\Q/\Z'(n)[-1]\quad\text{ for}\; n<0, $$
where $\Q/\Z'(n)=\colim_{p\not|m}\mu_m^{\otimes n}$ is the prime to 
$p$-part of $\Q/\Z(n)$. Then for any
$n$, there is are quasi-isomorphisms of complexes of \'etale
sheaves \cite[Prop. 6.7]{susvoe} and \cite{marcI},
$$ \Z/m(n) \cong
\begin{cases}
\mu_m^{\otimes n}[0] & p\not| m \\
\nu_{r}^n[-n] & m=p^r,
\end{cases}$$
where $\nu^n_{r}=W_r\Omega^n_{X,\log}$ is the logarithmic de
Rham-Witt sheaf. We let $H^i_{\cal M}(X,A(n))$ be the Nisnevich and 
$H^i_\et(X,A(n))$ be the \'etale
hypercohomology of $A(n)$, and abbreviate the Weil-\'etale
hypercohomology $H^i_W(X,\gamma^*A(n))$ by
$H^i_W(X,A(n))$. Then $H^i_\et(X,\Q(n))\cong H^i_\M(X,\Q(n))$. 
Theorems \ref{maincup}, \ref{kahnequal}, Corollary
\ref{rational}, and Proposition \ref{cupe} specialized to this
situation give:

\begin{theorem}
\label{ttt}
Let $X$ be a smooth variety over $\F_q$.

a) There is a long exact sequence
\begin{equation}
\label{mainsequence} \cdots \to H^i_\et(X,\Z(n))\to H^i_W(X,\Z(n))
\to H^{i-1}_\et(X,\Q(n)) @>\delta >> H^{i+1}_\et(X,\Z(n))\to
\cdots,
\end{equation}
where the map $\delta$ is the composition
$$ H^{i-1}_\et(X,\Q(n))\to
H^{i-1}_\et(X,\Q/\Z(n)) @>\cdot \epsilon
>> H^{i}_\et(X,\Q/\Z(n))@>\beta >> H^{i+1}_\et(X,\Z(n)). $$

b) For torsion coefficients, we have
$$ H^i_\et(X,\Z/m(n))\cong H^i_W(X,\Z/m(n)).$$

c) With rational coefficients,
$$H^i_W(X,\Q(n))\cong H^i_\M(X,\Q(n))\oplus
H^{i-1}_\M(X,\Q(n)).$$ The cup product with $e\in \Ext^1_G(\Z,\Z)$
is multiplication by the matrix
${0\,0\choose1\, 0}$.
\end{theorem}

\proofend

\begin{lemma}
\label{666} Let $\bar X$ be a smooth variety of dimension $d$ over
an algebraically closed field of characteristic $p$, and let
$l\not=p$.

a) If $n\geq d$, then $H^i_\et(\bar X,\Z(n))\cong H^i_\M(\bar
X,\Z(n))$, and the latter group is zero for $i>n+d$.

b) If $n<d$, then $H^i_\et(\bar X,\Z(n))$ is zero for $i>2d+1$,
torsion for $i>2n$, $p$-torsion free for $i>n+d+1$,
$p$-divisible for $i>n+d$ and $l$-divisible for $i>2d$.
\end{lemma}

\proof a) We consider the statement rationally, with prime to 
$p$-coefficients and $p$-power coefficients separately.
Let $\epsilon:\bar X_{\et}\to \bar X_{Zar}$ be the
canonical morphism of sites. Rationally, $R\epsilon_*\Q(n)_\et\cong
\Q(n)_\Zar$ for any $n$. For $p\not|m$,
it follows from Suslin \cite{susetale} that 
$R\epsilon_*\Z/m(n)_\et \cong \Z/m(n)_\Zar$ for $n\geq d$. 
Finally, $R\epsilon_*\Z/p^r(n)_\et\cong \Z/p^r(n)_\Zar$, 
because both sides are zero for $n>d$ by \cite{marcI}, and 
$R^j\epsilon_*\nu^d_r=0$  for $j>0$ by Gros-Suwa 
\cite[III Lemme 3.16]{grossuwa}. 
Hence we have $H^{i}_\et(\bar X,\Z(n))\cong H^i_\M(\bar X,\Z(n))$ 
for any $i$ and $n\geq d$. The latter group is zero for $i>d+n$ by 
definition.

b) Note first that $H^i_\et(\bar X,\Q(n))=H^i_\M(\bar X,\Q(n))=0$
for $i>2n$. For mod $p$ coefficients, $\Z/p(n)_\et \cong
\nu^n_1[-n]$ implies that 
$H^i_\et(\bar X,\Z/p(n))=H^{i-n}_\et(\bar X,\nu^n_1)=0$ for $i>d+n$, because
$\cd_p \bar X=d$. For mod $l$ coefficients,
$\Z/l(n)_\et\cong\mu^{\otimes n}_l[0]$ implies that $H^i_\et(\bar
X,\Z/l(n))=0$ for $i>2d$ because $\cd_l\bar X=2d$. The statement
now follows using the short exact sequence
$$ 0\to H^i_\et(\bar X,\Z(n))/l \to H^i_\et(\bar X,\Z/l(n)) \to
{}_l H^{i+1}_\et(\bar X,\Z(n)) \to 0.$$ \proofend

\begin{theorem}\label{bound}
Let $X$ be a smooth variety over a finite field of dimension $d$. Then
$H^i_W(X,\Z(n))=0$ for $i>\max\{2d+1,n+d+1\}$.
\end{theorem}

\proof In view of \eqref{descent} and the previous lemma, the
result is clear for $n\geq d$. Similarly, for $n<d$ and $i>2d+2$,
the group in question vanishes, because $H^i_\et(\bar X,\Z(n))=0$
for $i>2d+1$. It remains to show $H^{2d+2}_W(X,\Z(n))\cong
H^{2d+1}_\et(\bar X,\Z(n))_G=0$. By the lemma, $H^{2d+1}_\et(\bar
X,\Z(n))$ is a divisible torsion group, and this
property is inherited by its quotient group $H^{2d+2}_W(X,\Z(n))$.
On the other hand, by Theorem \ref{ttt} b), the limit of the
surjections $H^{2d+1}_W(X,\Z/l^r(n))\to
{}_{l^r}H^{2d+2}_W(X,\Z(n))$ gives a surjection $
H^{2d+1}_{cont}(X,\Z_l(n)) \to T_l H^{2d+2}_W(X,\Z(n))$ for every
prime number $l$. By Deligne's proof of the Weil conjectures,
$H^i_{cont}(X,\Z_l(n))$ is torsion for $i>n+d+1$, hence 
$T_l H^{2d+2}_W(X,\Z(n))=0$. But a divisible $l$-torsion group $A$ with
$T_lA=0$ is trivial. Indeed, any non-zero element $a\in A$
gives by divisibility rise to a non-zero element $(a_i)_i$ in the
Tate-module.
\proofend

We can give more explicit formulas for the case $n=0,1$:

\begin{prop}\label{small} Let $X$ be smooth and connected.

a) We have 
$$H^i_W(X,\Z)\cong
\begin{cases}
\Z &i=0,1;\\
H^{i-1}_\et(X,\Q/\Z) &i>2.
\end{cases}$$

b) If $X$ is proper, then $H^i_W(X,\Z)$ is finite for $i>1$.

c) We have 
$$H^i_W(X,\Z(1))=
\begin{cases}
0 &i=0;\\
{\cal O}(X)^\times &i=1;\\
H^{i-1}_\et(X,\Q/\Z(1))&i\geq 5,
\end{cases}$$
and there is an exact sequence
\begin{multline*}
0 \to \Pic X \to H^2_W(X,\Z(1))\to {\cal O}(X)^\times\otimes\Q \to
\Br X \to H^3_W(X,\Z(1)) \to\\
\NS X\otimes\Q \to H^3_\et(X,\Q/\Z(1)) \to H^4_W(X,\Z(1)) \to 0.
\end{multline*}

d) If $X$ is proper, then $H^2_W(X,\Z(1))\cong \Pic X$, the groups
$H^i_W(X,\Z(1))$ are finite for $i\geq 5$, and there is an
exact sequence $$0\to \Br X \to H^3_W(X,\Z(1)) \to\\
\NS X\otimes\Q \to H^3_\et(X,\Q/\Z(1)) \to H^4_W(X,\Z(1)) \to 0.$$

If $H^3_W(X,\Z(1))$ is finitely generated, then $Br(X)$ is finite,
we have a short exact sequence
$$0\to \Br X\to H^3_W(X,\Z(1))\to \NS X \to 0,$$
and $H^4_W(X,\Z(1))\cong H^3_\et(X,\Q/\Z(1))_{cotor}$ is finite.
\end{prop}

\proof a) Since $H^1_\et(\bar X,\Z)=0$, we have
$H^1_W(X,\Z)=H^0_\et(\bar X,\Z)_G=\Z$. From \eqref{mainsequence}
and $H^i_\et(X,\Q)=0$ for $i>0$ we get
$H^{i-1}_\et(X,\Q/\Z)@>\beta>\sim >H^i_\et(X,\Z)@>\sim>>H^i_W(X,\Z)$.

b) This follows from Deligne's proof of the Weil conjectures, and
Gabber's finiteness result \cite{gabber}.

c) Since $\Z(1)\cong {\Bbb G}_m[-1]$ and $H^i_\et(X,\Q(1))=0$ for 
$i\geq 3$, this follows from sequence \eqref{mainsequence} together
with $H^2_\et(X,\Z(1))=\Pic X$ and $H^3_\et(X,\Z(1))=\Br X$.

d) The exact sequence follows because ${\cal O}(X)^\times$ is
torsion, and the finiteness again follows from Deligne's and
Gabber's results. If $H^3_W(X,\Z(1))$ is finitely generated, then
its image in $\NS X \otimes \Q$ is a lattice isomorphic to the 
torsion free finitely generated group $\NS X$ in order for the 
quotient to be torsion. Finally,
\begin{multline*}
\corank H^3_\et(X,\Q_l/\Z_l(1))=\dim H^3_{cont}(X,\Q_l(1)) =\dim
H^2_{cont}(\bar X,\Q_l(1))_{\hat G} \\
=\dim H^2_{cont}(\bar X,\Q_l(1))^{\hat G} =\dim
H^2_{cont}(X,\Q_l(1))=\rk \NS X.
\end{multline*}
\proofend

\example By Artin-Schreier theory,
$H^1_\et({\Bbb A}^1_{\F_q},\Z/p) \subseteq H^1_\et({\Bbb
A}^1_{\F_q},\Q_p/\Z_p)\cong H^2_\et({\Bbb A}^1_{\F_q},\Z)$
contains an infinite number of copies of $\Z/p$. From this and the
exact sequence
$$ \Q\cong H^0_\et({\Bbb A}^1_{\F_q},\Q)\to H^2_\et({\Bbb
A}^1_{\F_q},\Z)\to H^2_W({\Bbb A}^1_{\F_q},\Z)\to 0$$ we see that
$H^2_W({\Bbb A}^1_{\F_q},\Z)$ is not finitely generated.

\begin{theorem}
\label{degree2d} Let $X$ be a connected smooth projective variety
of dimension $d$ over $\F_q$ with Albanese variety $\Alb$. Then
there is a commutative diagram
$$\begin{CD}
H^{2d}_W(X,\Z(d)) @>>> H^{2d}_\et(\bar X,\Z(d))^G
@>(deg',Alb)>\sim> \Z \oplus \Alb(\F_q)\\
@V e VV @VVV @VVV \\
H^{2d+1}_W(X,\Z(d)) @<\sim <<  H^{2d}_\et(\bar X,\Z(d))_G
@>deg>\sim >  \Z
\end{CD}$$
\end{theorem}

\proof 
By Lemma \ref{666},  $H^{2d}_\et(\bar X,\Z(d))\cong
\CH^d(\bar X)$. The kernel of the degree map 
$\CH^d(\bar X)@>deg >>\Z$ is divisible, and by Rojtman's theorem 
\cite{roitman, milneroitman} it agrees with the
$\bar \F_q$-rational torsion points $\Alb(\bar \F_q)$.
Since $\Alb(\bar \F_q)^G= \Alb(\F_q)$ is finite, we get 
$\Alb(\bar \F_q)_G=0$. The statement now follows from the short exact
sequence \eqref{descent} and Lemma \ref{multss} b).
\proofend

\section{Comparison to $l$-adic cohomology}
Fix a smooth projective variety $X$ over $\F_q$ and an
integer $n$. There are two fundamental conjectures on Weil motivic
cohomology. The first one is due to Lichtenbaum:

\begin{conj} {\bf L(X,n)}
\label{lichtconj}
For every $i$, the group $H^i_W(X,\Z(n))$ is finitely generated.
\end{conj}

Note that $H^i_W(X,\Z(n))$ may be not finitely generated if $X$
is not smooth or not projective. 
The homomorphisms $\Z(n)\otimes\Z_l \to \Z/l^r(n)$ in the derived
category of sheaves  of $\T_G$ are compatible, and induce a
morphism $\Z(n)\otimes\Z_l \to R\lim \Z/l^r(n)$, hence in view of
Lemma \ref{lim} upon applying
$R\gamma_*$ a map $c:R\gamma_*\Z(n)\otimes \Z_l @>>> R\lim
\Z/l^r(n)$ in the derived category of sheaves of $\T_{\hat G}$. In
view of \cite[Lemma 3.8]{kahngl} and Theorem \ref{kahnequal}, the
following conjecture is equivalent to Conjecture 3.2 of Kahn
\cite{kahngl}.

\begin{conj} {\bf K(X,n)}
\label{mainconj} For every prime $l$, and any $i$, the map $c$
induces an isomorphism
$$H^i_W(X,\Z(n))\otimes\Z_l@>\sim >> H^i_{cont}(X,\Z_l(n)).$$
\end{conj}

This implies in particular that Weil motivic cohomology is an
integral model for $l$-adic and
$p$-adic cohomology. Finally, there is the classical conjecture,
due to Tate (part 1,2 for $l\not=p$) and Beilinson (part 3):

\begin{conj} {\bf T(X,n)}
\label{tateconj}
For every prime $l$,
\begin{enumerate}
\item The cycle map $CH^n(X)\otimes \Q_l \to
H^{2n}_{cont}(\bar X,\Q_l(n))^{\hat G}$
is surjective.
\item The $\hat G$-module $H^{2n}_{cont}(\bar X,\Q_l(n))$ is semi-simple
at the eigenvalue $1$.
\item Rational and numerical equivalence agree with rational
coefficients.\end{enumerate}
\end{conj}

\begin{theorem}
\label{thmthm}\label{kahnthm} 
Let $X$ be a smooth projective
variety over $\F_q$, and $n$ an integer. Then
$$ K(X,n)+K(X,d-n) \Rightarrow L(X,n)\Rightarrow K(X,n)
\Rightarrow T(X,n).$$
Conversely, if $T(X,n)$ holds for {\rm all}
smooth and projective varieties $X$ over $\F_q$ and {\rm all}
$n$, then $K(X,n)$ holds for all $X$ and $n$.
\end{theorem}

\proof 
$K(X,n)+ K(X,d-n)\Rightarrow L(X,n)$: This has been
proved in by Kahn \cite{kahnneu}. 

\noindent $L(X,n)\Rightarrow K(X,n)$: 
By finite generation and Theorem \ref{ttt} b), we have
\begin{multline*}
$$H^i_W(X,\Z(n))\otimes \Z_l\cong \lim H^i_W(X,\Z(n))/l^r
\cong \lim H^i_W(X,\Z/l^r(n))\\
\cong \lim H^i_\et(X,\Z/l^r(n))
\cong H^i_{cont}(X,\Z_l(n)).
\end{multline*}

\noindent $K(X,n)\Rightarrow T(X,n)$: This is proved exactly as in 
\cite[Prop. 3.9]{kahngl}.

\noindent $T(-,-)\Rightarrow K(X,n):$ The hypothesis implies that
$H^i_\M(X,\Q(n))\cong H^i_\et(X,\Q(n))=0$ for $i\not=2n$, \cite{ich}.
One can now use the argument of \cite[Prop. 3.9]{kahngl}.
\proofend

\noindent{\bf Remark.}
1) Conjecture $K(X,n)$  for all smooth and proper varieties
implies the same statement for all smooth varieties as long as
$l\not=p$. This follows by using localization sequences for
$l$-adic cohomology and Weil-cohomology, and de Jong's theorem on
alterations, see \cite[Lemma 4.1]{alter} or \cite[Section 5]{kahn}
for details.


2) For $l\not=p$, an equivalent formulation of Conjecture
\ref{mainconj} is that the cycle map induces a quasi-isomorphism
$$R\gamma_* \Z(n)\otimes\Z_l @>\sim >> R\lim_{\hat G} \Z/l^r(n)$$
in the derived category of \'etale sheaves on smooth schemes 
over $\F_q$. The right hand term is quasi-isomorphic
to $R\gamma_* R\lim \Z/l^r(n)$ by Lemma \ref{lim}, but the above
quasi-isomorphism is not induced by a quasi-isomorphism between
$\Z(n)\otimes\Z_l$ and $R\lim \Z/l^r (n)$ in the derived 
category of Weil-\'etale sheaves on smooth schemes over $\F_q$. 
For example, ${\cal H}^2(\Z(1)\otimes\Z_l)=0$ by Hilbert's Theorem 90,
but one can show that $R^2\lim \Z/l^r(1)\not=0$.


3)  In view of the spectral sequences \eqref{ss1} and
\eqref{moress}, Conjecture $K(X,n)$ would follow if for all $i$
and $l$, the map
$H^i_{\hat G}(\bar X,\Z(n))\otimes\Z_l@>>> H^i_{\hat G}(\bar
X,(\Z/l^\cdot(n)))$ is an isomorphism of $\hat G$-modules. These
maps fit into a commutative diagram of short exact sequences of
$\hat G$-modules
$$\begin{CD}
R^1\gamma_*H^{i-1}_\et(\bar X,\Z(n))\otimes\Z_l
@>>> H^i_{\hat G}(\bar X,\Z(n))\otimes\Z_l
@>>> \gamma_*H^{i}_\et(\bar X,\Z(n))\otimes\Z_l \\
@VVV @VVV @VVV \\
R^1\gamma_*H^{i-1}_{cont}(\bar X,\Z_l(n)) @>>> H^i_{\hat G}(\bar
X,(\Z/l^\cdot(n)))@>>> \gamma_*H^{i}_{cont}(\bar X,\Z_l(n)).
\end{CD}$$
Thus $K(X,n)$ follows from isomorphisms of
$\hat G$-modules
\begin{align*}
H^{i}_\et(\bar X,\Z(n))\otimes\Z_l\cong \gamma_*H^{i}_\et(\bar
X,\Z(n))\otimes\Z_l & @>\sim >>
\gamma_*H^{i}_{cont}(\bar X,\Z_l(n))\\
H^{i-1}_\et(\bar X,\Z(n))\otimes\Q_l\cong
R^1\gamma_*H^{i-1}_\et(\bar X,\Z(n))\otimes\Z_l&@> \sim >>
R^1\gamma_*H^{i}_{cont}(\bar X,\Z_l(n)).
\end{align*}
The former isomorphism for $i=2n$ was Tate's original formulation of his
conjecture \cite{tatepoles}.

\section{Values of zeta-functions, Examples}
We can reformulate results of Milne \cite{values} to find expressions
for values of zeta functions as conjectured by Lichtenbaum
\cite{lalt}. Since
$e\in \Ext_{G}^1(\Z,\Z)$ satisfies $e^2=0$, the Weil \'etale
cohomology groups $H^*_W(X,\Z(n))$ form a complex under cup product
with $e$. By Theorems \ref{ttt} c) and \ref{bound},
the cohomology groups of this complex are torsion, and only finitely 
many are non-zero. For a complex $C^\cdot$ of abelian groups with finitely
many finite cohomology groups, one defines
$$ \chi(C^\cdot) := \prod_i |H^i(C^\cdot)|^{(-1)^i}.$$
Let $X$ be a smooth projective scheme over $\F_q$, and
$\zeta(X,s)=Z(X,q^{-s})$ be its zeta function. Following Milne
\cite{values}, we let
$$ \chi(X,{\cal O}_X,n)=\sum_{i\leq n,j\leq d}(-1)^{i+j}(n-i)
\dim H^j(X,\Omega^i).$$
The conclusion
of the following theorem has been proved by Lichtenbaum
for $n=0$ in \cite{licht}.

\begin{theorem}
\label{lvalues}
Let $X$ be a smooth projective variety such that $K(X,n)$
holds. Then the order $\rho_n$ of the pole of $\zeta(X,s)$
at $s=n$ is $\rk H^{2n}_W(X,\Z(n))$, and
\begin{equation}
\label{zeta}
\zeta(X,s)=\pm
(1-q^{n-s})^{-\rho_n}\cdot \chi(H^*_W(X,\Z(n)),e)\cdot
q^{\chi(X,{\cal O}_X,n)} \quad {\text as}\; s\to n.
\end{equation}
If $K(X,d-n)$ holds also, then
\begin{equation*}
\label{zeta+} \chi(H^*_W(X,\Z(n)),e)= \prod_i
|H^i_W(X,\Z(n))_{tor}|^{(-1)^i}\cdot R^{-1},
\end{equation*}
where $R$ is the determinant of the pairing
$$ H^{2n}_W(X,\Z(n)) \times H^{2(d-n)}_W(X,\Z(d-n))
\to H^{2d}_W(X,\Z(d)) @>deg(-\cup e)>> \Z .$$
\end{theorem}

\proof Since $K(X,n)$ implies semi-simplicity of $l$-adic
cohomology, the first formula follows by comparing to the formulas
for $l$-adic cohomology in \cite[Thm. 0.1]{values}.

Recall that for a short exact sequence $0\to A^\cdot \to B^\cdot
\to C^\cdot\to 0$ we have $\chi(A^\cdot)\cdot \chi(C^\cdot) =
\chi(B^\cdot)$. Thus it suffices to show that
$R^{-1}=\chi(H^*_W(X,\Z(n))/tor,e)$. By hypothesis, the latter complex
consists only of the upper map in the following commutative diagram
$$\begin{CD}
H^{2n}_W(X,\Z(n))/tor @>e>> H^{2n+1}_W(X,\Z(n))/tor \\
@Vd'VV @VdVV \\
\Hom(H^{2(d-n)}_W(X,\Z(d-n)),\Z)@=
\Hom(H^{2(d-n)}_W(X,\Z(d-n)),\Z).
\end{CD}$$
The maps $d'$ and $d$ are given by $d'(x)(y)=\deg(x\cdot y\cdot e)$ and
$d(x)(y)=\deg(x\cdot y)$, where $\deg$ is the map of
Theorem \ref{degree2d}. Comparing with
$l$-adic cohomology and using \cite[Lemma 5.3]{mmc}, one sees that
$d$ is an isomorphism. Hence
$$\chi(H^*_W(X,\Z(n))/tor,e)=\frac{|\ker e|}{|\coker e|}
=\frac{|\ker d'|}{|\coker d'|}=\frac{1}{R}.$$
\proofend

We give some explicit examples for varieties satisfying the
hypothesis of the previous theorem.

\begin{prop}
Conjecture $L(X,n)$ holds for $n\leq 0$. In particular,
\eqref{zeta} holds for all smooth projective $X$ and $n\leq 0$.
\end{prop}

\proof 
For $n=0$, this is Proposition \ref{small} b).
For $n<0$, the proposition follows because $H^i_{cont}(X,\Q_l(n))=0$ 
by Deligne's proof of the Weil conjectures, hence
$$ H^i_W(X,\Z(n))\otimes\Z_l := H^{i-1}_\et(X,\Q_l/\Z_l(n))\cong
H^i_{cont}(X,\Z_l(n))$$
is finitely generated and torsion.
\proofend

\begin{theorem}
Assume that $X$ is smooth and projective and that the cycle
map $\Pic X\otimes \Q_l \to H^2_{cont}(X,\Q_l(1))$
is surjective for some $l$. Then $K(X,1)$ holds.
In particular, \eqref{zeta} holds for $X$ and $n=1$.
\end{theorem}

\proof In view of Theorem \ref{ttt} b) we can verify $K(X,1)$ after
tensoring with
$\Q$. We have $H^i_\et(X,\Q(1))=0$ for $i\not=2$,
$H^2_\et(X,\Q(1))=\Pic X\otimes\Q$
and $H^i_{cont}(X,\Q_l(1))=0$ for
$i\not=2,3$. Hence we get the following diagram from Theorem \ref{ttt} c)
$$\begin{CD}
\Pic X \otimes\Q_l@>\sim >> H^{2}_W(X,\Z(1))\otimes\Q_l @>surj >>
H^{2}_{cont}(X,\Q_l(1))\\
@| @VeVV @VeV V  \\
\Pic X\otimes\Q_l @<\sim << H^{3}_W(X,\Z(1))\otimes\Q_l @>>>
H^{3}_{cont}(X,\Q_l(1)).
\end{CD}$$
In codimension $1$, rational and homological equivalence agree
rationally, hence the upper composition is an isomorphism. On the
other hand, by Milne \cite[Prop. 0.3]{values}, the surjectivity of
the cycle map for one $l$ implies semi-simplicity of
$H^2_{cont}(\bar X,\Q_l(1))$, for all $l$ including $l=p$, hence the right
vertical map is an isomorphism. 
\proofend

In particular, the conclusion holds for Hilbert modular
surfaces, Picard modular surfaces, Siegel modular threefolds,
and in characteristic at least $5$ for supersingular
and elliptic K3 surfaces \cite{tateseattle}.
We use Soul\'e's method to produce more examples.

\begin{prop}
Let $X=X_1\times\ldots \times X_d$ be a product of smooth
projective curves over $\F_q$, and let $n\leq 1$ or $n\geq d-1$.
Then $K(X,n)$ holds for $X$.
\end{prop}

\proof By \cite[Prop. 3.9]{kahngl}, it suffices to show that
$H^i_\M(X,\Z(n))\otimes\Q=0$ for $i<2n$, that
$H^{2n}_\M(X,\Z(n))\otimes\Q_l \cong H^{2n}_{cont}(X,\Q_l(n))$,
and that the Frobenius acts semi-simply at $1$ on
$H^{2n}_{cont}(\bar X,\Q_l(n))$. We essentially repeat the proof
of Soul\'e \cite[Thm. 3]{soule}, adapted to our situation.

Write $X_i\cong {\bf 1}\oplus X^+_i\oplus \L$ in the category of
Chow motives, where ${\bf 1}=\Spec \F_q$ and ${\Bbb P}^1\cong {\bf
1}\oplus \L$. Then $X$ is a direct sum of motives of the form $M =
\otimes_{s=1}^j X_{n_s}^+\otimes \L^k$, with $0\leq j+k\leq d$.
Such a motive $M$ has a Frobenius endomorphism $F_M$, and $F_M$
has a minimal polynomial $P_M(u)$ such that all roots of $P_M(u)$
have absolute value equal to $q^{\frac{j+2k}{2}}$ 
\cite[Prop. 3.1.2]{soule}. Since the Frobenius $F_M$ acts on
$H^i_\M(M,\Z(n))\otimes\Q$ as multiplication by $q^n$ 
\cite[Prop. 1.5.2]{soule}, we get $0=P_M(F_M)=P_M(q^n)$ on
$H^i_\M(M,\Z(n))\otimes\Q$. For $j+2k\not=2n$, $P_M(q^n)$ is
non-zero, hence $H^i_\M(M,\Z(n))\otimes\Q=0$. By our choice of $n$,
$j+2k=2n$ can only happen for $j=0$ or $j=2$.

If $j=0$, then $M=\L^n$, and by the projective bundle formula for
motivic cohomology,
$H^i_\M(\L^n,\Z(n))\otimes\Q_l=H^{i-2n}_\M(\F_q,\Z(0))\otimes\Q_l$.
The latter group is zero except for $i=2n$, in which case it is
isomorphic to $H^{2n}_{cont}(\bar\L^n,\Q_l(n))=\Q_l$. The Galois
group acts trivially, in particular semi-simply, on the latter
group.

If $j=2$, then $M=X^+\otimes Y^+\otimes \L^{n-1}$, and
$$H^i_\M(M,\Z(n))\otimes\Q_l=H^{i-2n+2}_\M(X^+\otimes Y^+,\Z(1))\otimes\Q_l
=H^{i-2n+1}_\M(X^+\otimes Y^+,{\Bbb G}_m)\otimes \Q_l.$$
The
latter group is zero for $i<2n$, because the group of global
sections of a projective variety over a finite field is finite. On
the other hand,
\begin{multline*}
H^{2n}_\M(M,\Z(n))\otimes\Q_l\cong \CH^1(X^+\otimes Y^+)\otimes \Q_l \\
\cong H^{2}_{cont}(X^+\otimes Y^+,\Q_l(1))\cong
H^{2n}_{cont}(M,\Q_l(n))
\end{multline*}
by Tate's theorem \cite{tate}. Tate's theorem also implies that
the Galois group acts semi-simply at $1$ on the module
$H^{2}_{cont}(\overline{X^+\otimes Y^+},\Q_l(1))$. For $l=p$, the same
statement follows by \cite[Prop. 0.3]{values}. 
\proofend

As in Soul\'e, let $A(k)$ be the subcategory of smooth projective
varieties generated by products of curves and the following
operations:

\begin{enumerate}
\item If $X$ and $Y$ are in $A(k)$, then $X\coprod Y$ is in $A(k)$.
\item If $Y$ is in $A(k)$, and there are morphisms
$c:X\to Y$ and $c':Y\to X$ in the category of Chow motives, such
that $c'\circ c:X\to X$ is multiplication by a constant, then
$X$ is in $A(k)$.
\item If  $k'$ is a finite extension of $k$,
and $X\times_kk'$ is in $A(k')$, then $X$ is in $A(k)$.
\item If $Y$ is a closed subscheme of $X$ and $Y$ and $X$ are
in $A(k)$, then the blow-up $X'$ of $X$ along $Y$ is in $A(k)$.
\end{enumerate}

\begin{theorem}
Let $X$ be a variety of dimension $d$ in $A(\F_q)$. Then $K(X,n)$
and $L(X,n)$ hold for $n\leq 1$ or $n\geq d-1$. In particular,
\eqref{zeta} and \eqref{zeta+} hold for $X$ and $n\leq 1$ or
$n\geq d-1$.
\end{theorem}

\proof The statement holds for products of curves, and it is clear
that if $X$ and $Y$ satisfy $K(X,n)$ then $X\coprod Y$ also does.
In 2) and 3), the map $H^i_W(X,\Z(n))\otimes\Q_l\to
H^i_{cont}(X,\Q_l(n))$ is a direct summand for the corresponding
map for $Y$ and $X\times_k k'$, respectively. Finally, if $X'$ is
the blow-up of $X$ along $Y$, and $Y$ has codimension $c$ in $X$,
then one has $X'= X\oplus (\oplus_{j=1}^{c-1}Y\otimes \L^j)$.
\proofend

In particular, the conclusion of theorem holds for
abelian varieties, unirational
varieties of dimension at most $3$, or Fermat hypersurfaces.

In \cite{kahnneu}, Kahn shows that $K(X,n)$ holds, if the Chow
motive of $X$ is in the subcategory of Chow motives generated by
abelian varieties and Artin motives, and if Tate's conjecture
holds for $X$. This applies in particular to products of elliptic
curves.

\end{document}